\documentclass[11pt,twoside]{article}
\usepackage[mathscr]{euscript}
\usepackage{secdot}
\usepackage[utf8]{inputenc}
\usepackage{amsmath, amsthm, amscd, amsfonts, amssymb, color}
\usepackage[bookmarksnumbered,plainpages]{hyperref}
\newcommand\mybox[1]{\parbox[t]{1\textwidth}{\raggedright$\displaystyle #1 $}}
\allowdisplaybreaks

\date{}

\begin{document}

\centerline{}

\centerline {\Large{\bf Common Fixed Point Theorem for Six Functions on}}
\centerline{}
\centerline {\Large{\bf Menger Probabilistic  Generalized Metric Space}}
\centerline{}
\centerline{\textbf{Sanjay Roy}}
\centerline{Department of Mathematics, Uluberia College}
\centerline{Uluberia, Howrah,  West Bengal, India}
\centerline{E-mail: sanjaypuremath@gmail.com}
\centerline{}

% My definition
\newcommand{\mvec}[1]{\mbox{\bfseries\itshape #1}}
\centerline{}
\newtheorem{Theorem}{\quad Theorem}[section]

\newtheorem{definition}[Theorem]{\quad Definition}

\newtheorem{theorem}[Theorem]{\quad Theorem}

\newtheorem{remark}[Theorem]{\quad Remark}

\newtheorem{corollary}[Theorem]{\quad Corollary}

\newtheorem{note}[Theorem]{\quad Note}

\newtheorem{lemma}[Theorem]{\quad Lemma}

\newtheorem{example}[Theorem]{\quad Example}

\newtheorem{notation}[Theorem]{\quad Notation}

\newtheorem{result}[Theorem]{\quad Result}
\newtheorem{conclusion}[Theorem]{\quad Conclusion}

\newtheorem{proposition}[Theorem]{\quad Proposition}
\newtheorem{prop}[Theorem]{\quad Property}

\begin{abstract}
\textbf{\emph{ The main aim of this paper is to find a unique common fixed point for six functions in a Menger probabilistic generalized metric space. For this purpose, we have defined the compatibility of three functions and established some required theorems.
}}
\end{abstract}

%=====================================
\section{Introduction}
%=====================================
The concept of probabilistic metric space was first introduced by K. Manger \cite{Menger} in 1942 in the name of statistical metric to find the distance between two points in  probabilistic approach. Then Schweizer and Sklar \cite{Schweizer1, Schweizer2, Schweizer3}  redefined this in another way. Then many works such as strong ideal convergence \cite{Sencimen},  convergence in PM space, completeness \cite{Sherwood}, Semigroups \cite{Elamrani} etc. had been done on probabilistic metric space.

The well-known fixed point theorem plays an important role to solve many problems such as the existence of solutions, the existence of orbits in dynamical systems, image processing, and economics etc.  So many authors \cite{Chang1, Hadzic12, Hadzic13, Hadzic14, Hicks1, Singh, Sharma} established the fixed point theorem for various types of contraction mappings on probabilistic metric spaces in the last few decades. Again many authors \cite{Chang2, Hadzic6, Hadzic11, Kutukcu, Mbarki} tried to find out the common fixed point of functions, satisfying certain conditions in the last few years and in 2013, R. Singh etl. \cite{Singh} established a sufficient condition for the existence of common fixed point of four functions.

 In 2014, C. Zhou etl.\cite{Zhou} have generalized the concept of probabilistic metric space and established a fixed point theorem in this space.
 
 Motivated by the works of R. Singh etl. \cite{Singh} and C. Zhou etl.\cite{Zhou},  we have explored the existence and uniqueness of a common fixed point for six self-mappings within the frame of a Menger probabilistic generalized metric space.
  
%=====================================
\section{Preliminaries}
%=====================================
\smallskip\hspace{.6 cm} In this section, we recall some definitions and results in Menger Probabilistic generalized metric space.
\begin{definition}\cite{Schweizer3}
A binary relation $\ast$ on $[0,1]$ is said to be a triangular norm or $t$-norm if the following condition are satisfied:\\
$(i)$  $x\ast y=y\ast x$,\\
$(ii)$  $x\ast(y\ast z)=(x\ast y)\ast z$,\\
$(iii)$  $x\ast y\leq x\ast z$ whenever $y\leq z$\\
$(iv)$  $x\ast 1= x$ for all $x, y, z\in [0, 1]$.
\end{definition}

\begin{definition}
A distribution function is a function $F:[-\infty,\, \infty]\rightarrow [0, 1]$ which is left continuous on $\mathbb{R}$, non-decreasing and $F(-\infty)=0$, $F(\infty)=1$.
\end{definition}

\begin{definition}
The Dirac distribution function $H_a:[-\infty,\,\infty]\rightarrow [0, 1]$ is defined for $a\in[-\infty,\, \infty)$ by\\
\begin{center}
$H_a(u) = 
\begin{cases} 
0 & \text{if } u \in [-\infty,\, a], \\
1 & \text{if } u\in(a,\,\infty], 
\end{cases}$
\end{center}
and for $a=\infty$ by\\
\begin{center}
$H_\infty(u) = 
\begin{cases} 
0 & \text{if } u \in [-\infty,\, \infty), \\
1 & \text{if } u=\infty. 
\end{cases}$
\end{center}
\end{definition}

\begin{definition}
A distance distribution function $F:[-\infty,\, \infty]\rightarrow [0, 1]$ is a distribution function with $F(0)=0$. The family of all distance distribution functions is denoted by $\bigtriangleup^+$. 
\end{definition}

\begin{definition}\cite{Zhou}\label{D1}
A Menger probabilistic generalized metric space or Menger probabilistic G-metric space (briefly, a PGM space)  is a triple $(X, \mathcal{G}, \ast)$ where $X$ is a non empty set, $\ast$ is a continuous $t$-norm and $\mathcal{G}: X\times X\times X\rightarrow \bigtriangleup^+$, given by $(x, y, z)\mapsto G_{x,\, y,\, z}$ such that the following conditions are satisfied for all $x,y,z$ in $X$:\\
$(i)$ $G_{x,\, y,\, z}(t)= 1$ for all $x, y, z\in X$ and $t>0$ if and only if $x=y=z$,\\
$(ii)$ $G_{x,\, x,\, y}(t)\geq G_{x,\, y,\, z}(t)$ for all $x, y\in X$ with $z\neq y$ and $t> 0$,\\
$(iii)$ $G_{x,\, y,\, z}(t)=G_{x,\, z,\, y}(t)=G_{y,\, x,\, z}(t)=\cdots$ $($ symmetry in all three variables $)$,\\
$(iv)$ $G_{x,\, y,\, z}(t+s)\geq G_{x,\, a,\, a}(t)\ast G_{a,\, y,\, z}(s)$ for all $x, y, z, a\in X$ and $s, t\geq 0$.
\end{definition}

\begin{definition}\cite{Zhou}
Let $(X, \mathcal{G}, \ast)$ be a PGM- space and $x_0\in X$. For any $\epsilon>0$ and $\delta$ with $0<\delta<1$, an $(\epsilon,\,\delta)$- neighbourhood of $x_0$ is denoted by $N_{x_0}(\epsilon, \delta)$ and is defined by
\begin{center}
$N_{x_0}(\epsilon, \delta)=\{y\in X:\; G_{x_0,\, y,\, y}(\epsilon)>1-\delta$ and $G_{y,\, x_0,\, x_0}(\epsilon)>1-\delta\}$
\end{center}
\end{definition}

\begin{definition}\cite{Zhou}
A sequence $\{x_n\}$ in a PGM-space $(X, \mathcal{G},\ast)$ is said to be convergent to a point $x\in X$ if for any $\epsilon>0$ and $0<\delta<1$, there exists a positive integer $M_{\epsilon,\,\delta}$ such that $x_n\in N_x(\epsilon,\,\delta)$ whenever $n>M_{\epsilon,\,\delta}$.

A sequence $\{x_n\}$ in a PGM-space $(X, \mathcal{G},\ast)$ is said to be a Cauchy sequence if for any $\epsilon>0$ and $0<\delta<1$, there exists a positive integer $M_{\epsilon,\,\delta}$ such that $G_{x_n, x_m, x_l}(\epsilon)>1-\delta$ whenever $m, n, l>M_{\epsilon,\,\delta}$.

A PGM-space $(X, \mathcal{G}, \ast)$ is said to be complete if every Cauchy sequence in $X$ converges to a point in $X$.
\end{definition}

\begin{theorem}\cite{Zhou}\label{th3}
Let $(X, \mathcal{G}, \ast)$ be a PGM- space. Let $\{x_n\},\, \{y_n\},\, \{z_n\}$ be sequences in $X$ and $x, y, z\in X$. If $x_n\rightarrow x,\; y_n\rightarrow y,\; z_n\rightarrow z$ as $n\rightarrow\infty$ then for any $t>0$, $G_{x_n,\;y_n,\;z_n}(t)\rightarrow G_{x,\; y,\; z}(t)$ as $n\rightarrow\infty$.
\end{theorem}

\begin{definition}
Let $(X, \mathcal{G}_1, \ast_1)$ and $(Y, \mathcal{G}_2, \ast_2)$ be two PGM- space. A function $f:(X, \mathcal{G}_1, \ast_1)\rightarrow (Y, \mathcal{G}_2, \ast_2)$ is said to be continuous at a point $z\in X$ if for any $\epsilon_2>0$ and $0<\delta_2<1$, there exist $\epsilon_1>0$ and $0<\delta_1<1$ such that $x\in N_z(\epsilon_1, \delta_1)$ implies that $fx\in N_{fz}(\epsilon_2, \delta_2)$.
\end{definition}

%=================================================
\section{Main Result}
%=================================================
\begin{definition}
Let $(X, \mathcal{G}, \ast)$ be a PGM- space and $A, B, C:X\rightarrow X$ be three self functions.
The triple $[A,B,C]$ is called compatible if for every sequence $\{x_n\}$ in $X$ with $\lim _{n\rightarrow\infty}Ax_n=\lim _{n\rightarrow\infty}Bx_n=\lim _{n\rightarrow\infty}Cx_n$ implies that $G_{\alpha,\; \beta,\;\beta}\rightarrow H_0$ for all $\alpha, \beta\in\{ABx_n,\; BCx_n,\;CAx_n\}$.

We say that the pair $[A,\, B]$ is compatible if the triple $[A,\, B,\, C]$ is compatible for $B=C$.
\end{definition}

\begin{remark}
If $ABx= BCx=CAx$ for all $x\in X$ then $[A, B, C]$ is compatible.
\end{remark}

\begin{theorem}\label{th1}
Let $(X, \mathcal{G}, \ast)$ be a PGM- space and $f:X\rightarrow X$ be continuous. Let $\{x_n\}$ be a sequence in $X$ such that $\lim _{n\rightarrow\infty}x_n=z$ for some $z\in X$. Then $\lim _{n\rightarrow\infty}fx_n=fz$ that is, $G_{fz,\,fx_n,\,fx_n}\rightarrow H_0$ and $G_{fx_n,\,fz,\,fz}\rightarrow H_0$.
\end{theorem}

\proof Since $f$ is continuous at $z$, for any $\epsilon>0$ and $0<\delta<1$ there exist $\epsilon_1>0$ and $0<\delta_1<1$ such that $x\in N_z(\epsilon_1, \delta_1)$ implies that $fx\in N_{fz}(\epsilon, \delta)$. Since $\lim _{n\rightarrow\infty}x_n=z$, there exists a natural number $m$ such that $x_n\in N_z(\epsilon_1, \delta_1)$ for all $n\geq m$.
So, if $n\geq m$ then $fx_n\in N_{fz}(\epsilon, \delta)$. Therefore $G_{fz,\,fx_n,\,fx_n}(\epsilon)>1-\delta$ and $G_{fx_n,\,fz,\,fz}(\epsilon)>1-\delta$ for all $n\geq m$. Hence $G_{fz,\,fx_n,\,fx_n}\rightarrow H_0$ and $G_{fx_n,\,fz,\,fz}\rightarrow H_0$. 

\begin{theorem}\label{th2}
Let $(X, \mathcal{G}, \ast)$ be a PGM- space and $A, B, C:X\rightarrow X$ where $B$ and $C$ are continuous. Let $\{x_n\}$ be a sequence in $X$ such that $\lim _{n\rightarrow\infty}Ax_n=\lim _{n\rightarrow\infty}Bx_n=\lim _{n\rightarrow\infty}Cx_n=z$ for some $z\in X$. If $[A,B,C]$ is compatible then $ABx_n\rightarrow Bz$ and $ACx_n\rightarrow Cz$.
\end{theorem}

\proof Let $\epsilon>0$ and $0<\delta<1$. Since $(X, \mathcal{G}, \ast)$ is a PGM- space,\\
$G_{Bz,\; ABx_n,\; ABx_n}(\epsilon)\geq G_{Bz,\; BCx_n,\; BCx_n}(\frac{\epsilon}{2})\ast G_{BCx_n,\; ABx_n, \;ABx_n}(\frac{\epsilon}{2})$ and\\
 $G_{ABx_n,\; Bz,\; Bz}(\epsilon)\geq G_{ABx_n,\; BCx_n,\; BCx_n}(\frac{\epsilon}{2})\ast G_{BCx_n,\; Bz,\; Bz}(\frac{\epsilon}{2})$.\\
Since $\lim _{n\rightarrow\infty}Ax_n=\lim _{n\rightarrow\infty}Bx_n=\lim _{n\rightarrow\infty}Cx_n=z$ and $[A,B,C]$ is compatible, $G_{\alpha,\; \beta,\;\beta}\rightarrow H_0$ for all $\alpha, \beta\in\{ABx_n,\; BCx_n,\;CAx_n\}$. \\
So, $G_{BCx_n,\; ABx_n, \;ABx_n}\rightarrow H_0$
 and $G_{ABx_n,\; BCx_n,\; BCx_n}\rightarrow H_0$. $\hfill\cdots(i)$ \\
  Again, since $\lim _{n\rightarrow\infty}Cx_n=z$ and $B$ is continuous, $\lim _{n\rightarrow\infty}BCx_n=Bz$ $[$ by the Theorem \ref{th1} $]$.\\
So, $G_{Bz,\; BCx_n,\; BCx_n}\rightarrow H_0$
and $G_{BCx_n,\; Bz,\; Bz}\rightarrow H_0$. $\hfill\cdots(ii)$\\
 Then from $(i)$ and $(ii)$, we can say that there exists $m\in\mathbb{N}$  such that\\  $G_{BCx_n,\; ABx_n, \;ABx_n}(\frac{\epsilon}{2})>1-\delta$, $G_{ABx_n,\; BCx_n,\; BCx_n}(\frac{\epsilon}{2})>1-\delta$, $G_{Bz,\; BCx_n,\; BCx_n}(\frac{\epsilon}{2})>1-\delta$ and $G_{BCx_n,\; Bz,\; Bz}(\frac{\epsilon}{2})>1-\delta$ for all $n\geq m$. Thus $G_{Bz,\; ABx_n,\; ABx_n}(\epsilon)\geq 1-\delta$ and $G_{ABx_n,\; Bz,\; Bz}(\epsilon)\geq 1-\delta$ for all $n\geq m$. Hence $ABx_n\rightarrow Bz$. Similarly, we can show that $ACx_n\rightarrow Cz$.

\begin{theorem}
Let $(X, \mathcal{G}, \ast)$ be a complete PGM- space  and $a*a\geq a$ for all $a\in [0, 1]$. Let $A, B, C, D, S$ and $T$ be self mappings on $X$ such that the following conditions are satisfied:\\
$(a)$ there exists a sequence $\{x_n\}_{n=0}^\infty$ of distinct elements of $X$ such that $Tx_{3n+1}=Ax_{3n}$, $Dx_{3n+2}=Bx_{3n+1}$, $Sx_{3n+3}=Cx_{3n+2}$, $Ax_{3n}\neq Bx_{3n+1}$, $Bx_{3n+1}\neq Cx_{3n+2}$ and $Cx_{3n+2}\neq Ax_{3n+3}$ for all $n= 0, 1, 2, \cdots$\\
$(b)$ $[A, S]$, $[B, T]$ and $[C, D]$ are compatible.\\
$(c)$ $D$, $S$ and $T$ are continuous.\\
$(d)$ There exists a number $k\in(0, \frac{1}{2}]$ such that\\ 
$G_{Ax,\,By,\, Cz}(kt)\geq G_{Sx,\, Ty,\, Dz}(t)*G_{Sx,\, Ax,\, Dz}(t)*G_{Ax,\, By,\, Cz}(t)*G_{Ty,\, By,\, Cz}(t)*G_{Sx,\, Cz,\, Dz}(2t)$ for all $x, y, z\in X$ and $t>0$.\\
Then $A, B, C, D, S, T$ have a unique common fixed point.
\end{theorem}

\proof  We first construct a sequence $\{y_n\}$ as follows:\\
$y_{3n}=Tx_{3n+1}=Ax_{3n}$,\\ $y_{3n+1}= Dx_{3n+2}=Bx_{3n+1}$ and\\ $y_{3n+2}=Sx_{3n+3}=Cx_{3n+2}$ for $n=0, 1, 2, \cdots$.\\
So, by our assumption $(a)$, we can say that  $\{y_n\}_{n=0}^\infty$ is a sequence in $X$ such that $y_n\neq y_{n+1}$ for all $n\in\mathbb{N}$ .
Now for $t>0$ we get
\begin{align*}
&G_{y_{3n},\, y_{3n+1},\, y_{3n+2}}(kt)\\
&= G_{Ax_{3n},\, Bx_{3n+1},\, Cx_{3n+2}}(kt)\\
&\geq \mybox{G_{Sx_{3n},\, Tx_{3n+1},\, Dx_{3n+2}}(t)*G_{Sx_{3n},\, Ax_{3n},\, Dx_{3n+2}}(t)*G_{Ax_{3n},\, Bx_{3n+1},\, Cx_{3n+2}}(t)*G_{Tx_{3n+1},\, Bx_{3n+1},\, Cx_{3n+2}}(t)*G_{Sx_{3n},\, Cx_{3n+2},\, Dx_{3n+2}}(2t)}\\
&\geq \mybox{G_{y_{3n-1},\, y_{3n},\, y_{3n+1}}(t)*G_{y_{3n-1},\, y_{3n},\, y_{3n+1}}(t)*G_{y_{3n},\, y_{3n+1},\, y_{3n+2}}(t)* G_{y_{3n},\, y_{3n+1},\, y_{3n+2}}(t)*G_{y_{3n-1},\, y_{3n+2},\, y_{3n+1}}(2t)}\\
&\geq G_{y_{3n-1},\, y_{3n},\, y_{3n+1}}(t)*G_{y_{3n},\, y_{3n+1},\, y_{3n+2}}(t)*
G_{y_{3n+2},\, y_{3n+1},\, y_{3n-1}}(2t)\;\; [\text{ as } a* a\geq a]\\
&\geq G_{y_{3n-1},\, y_{3n},\, y_{3n+1}}(t)*G_{y_{3n},\, y_{3n+1},\, y_{3n+2}}(t)*
 G_{y_{3n+2},\, y_{3n},\, y_{3n}}(t)*
  G_{y_{3n},\, y_{3n+1},\, y_{3n-1}}(t)\\
&\geq \mybox{G_{y_{3n-1},\, y_{3n},\, y_{3n+1}}(t)*G_{y_{3n},\, y_{3n+1},\, y_{3n+2}}(t)* G_{y_{3n},\, y_{3n+1},\, y_{3n+2}}(t)*G_{y_{3n},\, y_{3n+1},\, y_{3n-1}}(t)\;[\text{ as } y_{3n+1}\neq y_{3n+2}\; ]}\\
&\geq G_{y_{3n-1},\, y_{3n},\, y_{3n+1}}(t)*G_{y_{3n},\, y_{3n+1},\, y_{3n+2}}(t)* G_{y_{3n-1},\, y_{3n},\, y_{3n+1}}(t) \\
 &\geq \mybox{G_{y_{3n-1},\, y_{3n},\, y_{3n+1}}(t)*G_{y_{3n},\, y_{3n+1},\, y_{3n+2}}(t)\text{ for all } n\in\mathbb{N} \hfill \cdots (i)}
\end{align*}

\begin{align*}
&G_{y_{3n+1},\, y_{3n+2},\, y_{3n+3}}(kt)\\
&=G_{y_{3n+3},\, y_{3n+1},\, y_{3n+2}}(kt)= G_{Ax_{3n+3},\, Bx_{3n+1},\, Cx_{3n+2}}(kt)\\
&\geq \mybox{G_{Sx_{3n+3},\, Tx_{3n+1},\, Dx_{3n+2}}(t)*G_{Sx_{3n+3},\, Ax_{3n+3},\, Dx_{3n+2}}(t)*G_{Ax_{3n+3},\, Bx_{3n+1},\, Cx_{3n+2}}(t)* G_{Tx_{3n+1},\, Bx_{3n+1},\, Cx_{3n+2}}(t)*G_{Sx_{3n+3},\, Cx_{3n+2},\, Dx_{3n+2}}(2t)}\\
&\geq \mybox{G_{y_{3n+2},\, y_{3n},\, y_{3n+1}}(t)*G_{y_{3n+2},\, y_{3n+3},\, y_{3n+1}}(t)*G_{y_{3n+3},\, y_{3n+1},\, y_{3n+2}}(t)*G_{y_{3n},\, y_{3n+1},\, y_{3n+2}}(t)*G_{y_{3n+2},\, y_{3n+2},\, y_{3n+1}}(2t)}\\
&\geq \mybox{G_{y_{3n},\, y_{3n+1},\, y_{3n+2}}(t)*G_{y_{3n+1},\, y_{3n+2},\, y_{3n+3}}(t)*G_{y_{3n+1},\, y_{3n+2},\, y_{3n+3}}(t) *G_{y_{3n+1},\, y_{3n+2},\, y_{3n+2}}(2t)}\\
&\geq \mybox{G_{y_{3n},\, y_{3n+1},\, y_{3n+2}}(t)*G_{y_{3n+1},\, y_{3n+2},\, y_{3n+3}}(t)* G_{y_{3n+1},\, y_{3n},\, y_{3n}}(t)*G_{y_{3n},\, y_{3n+2},\, y_{3n+2}}(t)}\\
&\geq \mybox{G_{y_{3n},\, y_{3n+1},\, y_{3n+2}}(t)*G_{y_{3n+1},\, y_{3n+2},\, y_{3n+3}}(t)* G_{y_{3n+1},\, y_{3n},\, y_{3n+2}}(t) *G_{y_{3n},\, y_{3n+2},\, y_{3n+1}}(t)\;\; [\text{ as } y_{3n+1}\neq y_{3n+2} \text{ and }  y_{3n}\neq y_{3n+1}\; ]}\\
&\geq \mybox{G_{y_{3n},\, y_{3n+1},\, y_{3n+2}}(t)*G_{y_{3n+1},\, y_{3n+2},\, y_{3n+3}}(t)* G_{y_{3n},\, y_{3n+1},\, y_{3n+2}}(t) *G_{y_{3n},\, y_{3n+1},\, y_{3n+2}}(t)}\\
&\geq \mybox{G_{y_{3n},\, y_{3n+1},\, y_{3n+2}}(t)*G_{y_{3n+1},\, y_{3n+2},\, y_{3n+3}}(t)\text{ for all } n\in\mathbb{N}\cup\{0\}\hfill \cdots (ii)}
\end{align*}
\begin{align*}
&G_{y_{3n+2},\, y_{3n+3},\, y_{3n+4}}(kt)\\
&=G_{y_{3n+3},\, y_{3n+4},\, y_{3n+2}}(kt)= G_{Ax_{3n+3},\, Bx_{3n+4},\, Cx_{3n+2}}(kt)\\
&\geq \mybox{G_{Sx_{3n+3},\, Tx_{3n+4},\, Dx_{3n+2}}(t)*G_{Sx_{3n+3},\, Ax_{3n+3},\, Dx_{3n+2}}(t)*G_{Ax_{3n+3},\, Bx_{3n+4},\, Cx_{3n+2}}(t)* G_{Tx_{3n+4},\, Bx_{3n+4},\, Cx_{3n+2}}(t)*G_{Sx_{3n+3},\, Cx_{3n+2},\, Dx_{3n+2}}(2t)}\\
&\geq \mybox{G_{y_{3n+2},\, y_{3n+3},\, y_{3n+1}}(t)*G_{y_{3n+2},\, y_{3n+3},\, y_{3n+1}}(t)*G_{y_{3n+3},\, y_{3n+4},\, y_{3n+2}}(t)* G_{y_{3n+3},\, y_{3n+4},\, y_{3n+2}}(t)*G_{y_{3n+2},\, y_{3n+2},\, y_{3n+1}}(2t)}\\
&= \mybox{G_{y_{3n+1},\, y_{3n+2},\, y_{3n+3}}(t)*G_{y_{3n+1},\, y_{3n+2},\, y_{3n+3}}(t)*G_{y_{3n+2},\, y_{3n+3},\, y_{3n+4}}(t)*G_{y_{3n+2},\, y_{3n+3},\, y_{3n+4}}(t) *G_{y_{3n+1},\, y_{3n+2},\, y_{3n+2}}(2t)}\\
&\geq \mybox{G_{y_{3n+1},\, y_{3n+2},\, y_{3n+3}}(t)*
G_{y_{3n+2},\, y_{3n+3},\, y_{3n+4}}(t) *
G_{y_{3n+1},\, y_{3n+2},\, y_{3n+2}}(2t)}\\
&\geq \mybox{G_{y_{3n+1},\, y_{3n+2},\, y_{3n+3}}(t)*G_{y_{3n+2},\, y_{3n+3},\, y_{3n+4}}(t) *G_{y_{3n+2},\, y_{3n+3},\, y_{3n+3}}(t)\ast G_{y_{3n+1},\, y_{3n+2},\, y_{3n+3}}(t) \text{ by the Definition }\ref{D1}}\\
&\geq \mybox{G_{y_{3n+1},\, y_{3n+2},\, y_{3n+3}}(t)*G_{y_{3n+2},\, y_{3n+3},\, y_{3n+4}}(t) *G_{y_{3n+1},\, y_{3n+2},\, y_{3n+3}}(t) * G_{y_{3n+1},\, y_{3n+2},\, y_{3n+3}}(t)\;\;[\text{ as } y_{3n+1}\neq y_{3n+2}\; ]}\\
 &\geq \mybox{G_{y_{3n+1},\, y_{3n+2},\, y_{3n+3}}(t)*G_{y_{3n+2},\, y_{3n+3},\, y_{3n+4}}(t) \;\;[ \text{ as } a*a\geq a\; ] \text{ for all } n\in\mathbb{N}\cup\{0\} \hfill \cdots (iii)}
\end{align*}
Therefore, from $(i),\, (ii)$ and $(iii)$   we see that for all $n\in \mathbb{N}\cup\{0\}$ and for all $t> 0$, \\ $G_{y_{n+1},\, y_{n+2},\, y_{n+3}}(kt)\geq G_{y_{n},\, y_{n+1},\, y_{n+2}}(t)*G_{y_{n+1},\, y_{n+2},\, y_{n+3}}(t) $\\ 
that is,
$G_{y_{n+1},\, y_{n+2},\, y_{n+3}}(t)\geq G_{y_{n},\, y_{n+1},\, y_{n+2}}(k^{-1}t)*G_{y_{n+1},\, y_{n+2},\, y_{n+3}}(k^{-1}t)$ $\hfill\cdots(iv)$ \\\\
So, from $(iv)$ we get, $G_{y_{n+1},\, y_{n+2},\, y_{n+3}}(k^{-1}t)\geq G_{y_{n},\, y_{n+1},\, y_{n+2}}(k^{-2}t)*G_{y_{n+1},\, y_{n+2},\, y_{n+3}}(k^{-2}t)$.\\
Therefore, 
\begin{align*} 
&G_{y_{n+1},\, y_{n+2},\, y_{n+3}}(t)\\
&\geq \mybox{G_{y_{n},\, y_{n+1},\, y_{n+2}}(k^{-1}t)* G_{y_{n},\, y_{n+1},\, y_{n+2}}(k^{-2}t)*G_{y_{n+1},\, y_{n+2},\, y_{n+3}}(k^{-2}t)}\\
&\geq G_{y_{n},\, y_{n+1},\, y_{n+2}}(k^{-1}t)*G_{y_{n},\, y_{n+1},\, y_{n+2}}(k^{-1}t)*G_{y_{n+1},\, y_{n+2},\, y_{n+3}}(k^{-2}t) \\
&\smallskip\hspace{1 cm}[\text{ since }k^{-2}> k^{-1}, G_{y_{n},\, y_{n+1},\, y_{n+2}}(k^{-2}t)\geq G_{y_{n},\, y_{n+1},\, y_{n+2}}(k^{-1}t)]\\
&\geq \mybox{G_{y_{n},\, y_{n+1},\, y_{n+2}}(k^{-1}t)*G_{y_{n+1},\, y_{n+2},\, y_{n+3}}(k^{-2}t).}
\end{align*}
Again, from $(iv)$ we can write, \\
$G_{y_{n+1},\, y_{n+2},\, y_{n+3}}(k^{-2}t)\geq G_{y_{n},\, y_{n+1},\, y_{n+2}}(k^{-3}t)*G_{y_{n+1},\, y_{n+2},\, y_{n+3}}(k^{-3}t).$ So,
\begin{align*}
&G_{y_{n+1},\, y_{n+2},\, y_{n+3}}(t)\\
&\geq \mybox{G_{y_{n},\, y_{n+1},\, y_{n+2}}(k^{-1}t)* G_{y_{n},\, y_{n+1},\, y_{n+2}}(k^{-3}t)*G_{y_{n+1},\, y_{n+2},\, y_{n+3}}(k^{-3}t)}\\
&\geq \mybox{G_{y_{n},\, y_{n+1},\, y_{n+2}}(k^{-1}t)*G_{y_{n},\, y_{n+1},\, y_{n+2}}(k^{-1}t)*G_{y_{n+1},\, y_{n+2},\, y_{n+3}}(k^{-3}t)}\\
&\smallskip\hspace{1 cm} [ \text{ since }k^{-3}> k^{-1}, G_{y_{n},\, y_{n+1},\, y_{n+2}}(k^{-3}t)\geq G_{y_{n},\, y_{n+1},\, y_{n+2}}(k^{-1}t)]\\
&\geq G_{y_{n},\, y_{n+1},\, y_{n+2}}(k^{-1}t)*G_{y_{n+1},\, y_{n+2},\, y_{n+3}}(k^{-3}t).
\end{align*}
Proceeding in this way, we get \\
$G_{y_{n+1},\, y_{n+2},\, y_{n+3}}(t)\geq G_{y_{n},\, y_{n+1},\, y_{n+2}}(k^{-1}t)*G_{y_{n+1},\, y_{n+2},\, y_{n+3}}(k^{-m}t)$ for all $m\in \mathbb{N}$ and $n\in \mathbb{N}\cup\{0\}$.
Since $G_{y_{n+1},\, y_{n+2},\, y_{n+3}}(k^{-m}t)\rightarrow 1$ as $m\rightarrow\infty$ and $*$ is continuous,\\ 
$G_{y_{n+1},\, y_{n+2},\, y_{n+3}}(t)\geq G_{y_{n},\, y_{n+1},\, y_{n+2}}(k^{-1}t)$ for all $t>0$ and $n\in \mathbb{N}\cup\{0\}$.\\
Replacing $n$ by $n-1$ and $t$ by $k^{-1}t$, we get from above
$G_{y_{n},\, y_{n+1},\, y_{n+2}}(k^{-1}t)\geq G_{y_{n-1},\, y_{n},\, y_{n+1}}(k^{-2}t)$.
Proceeding in this way, we get\\
$G_{y_{n+1},\, y_{n+2},\, y_{n+3}}(t)\geq G_{y_{n},\, y_{n+1},\, y_{n+2}}(k^{-1}t)\geq G_{y_{n-1},\, y_{n},\, y_{n+1}}(k^{-2}t)\geq\\
\smallskip\hspace{9cm}\cdots\geq G_{y_{0},\, y_{1},\, y_{2}}(k^{-(n+1)}t)$.\\
We now show that $\{y_n\}$ is a Cauchy sequence. Let $p, q\in \mathbb{N}$ with $p\leq q$. Then\\
\begin{align*}
&G_{y_n,\,y_{n+p},\,y_{n+q}}(t)\\
&\geq \mybox{G_{y_n,\,y_{n+1},\,y_{n+1}}(\frac{t}{2})\ast G_{y_{n+1},\,y_{n+p},\,y_{n+q}}(\frac{t}{2})}\\
&\geq \mybox{G_{y_n,\,y_{n+1},\,y_{n+1}}(\frac{t}{2})\ast G_{y_{n+1},\,y_{n+2},\,y_{n+2}}(\frac{t}{2^2})\ast G_{y_{n+2},\,y_{n+p},\,y_{n+q}}(\frac{t}{2^2})}\\
&\geq \mybox{G_{y_n,\,y_{n+1},\,y_{n+1}}(\frac{t}{2})\ast G_{y_{n+1},\,y_{n+2},\,y_{n+2}}(\frac{t}{2^2})\ast \cdots \ast G_{y_{n+p-1},\,y_{n+p},\,y_{n+p}}(\frac{t}{2^p})\ast G_{y_{n+p},\,y_{n+p},\,y_{n+q}}(\frac{t}{2^p}).\hfill\cdots(v)}
\end{align*}
Now
\begin{align*}
&G_{y_{n+p},\,y_{n+p},\,y_{n+q}}(\frac{t}{2^p})\\
&\geq \mybox{G_{y_{n+p},\,y_{n+q-1},\,y_{n+q}}(\frac{t}{2^p})\;\; [\; \text{ as } y_{n+q-1}\neq y_{n+q} \; ]}\\
&\geq \mybox{G_{y_{n+p},\,y_{n+p+1},\,y_{n+p+1}}(\frac{t}{2^{p+1}})\ast G_{y_{n+p+1},\,y_{n+q-1},\,y_{n+p}}(\frac{t}{2^{p+1}})}\\
&\geq \mybox{G_{y_{n+p},\,y_{n+p+1},\,y_{n+p+1}}(\frac{t}{2^{p+1}})\ast G_{y_{n+p+1},\,y_{n+p+2},\,y_{n+p+2}}(\frac{t}{2^{p+2}})\ast\cdots\ast G_{y_{n+q-2},\,y_{n+q-1},\,y_{n+q-1}}(\frac{t}{2^{q-1}}) \ast G_{y_{n+q-1},\,y_{n+q-1},\,y_{n+q}}(\frac{t}{2^{q-1}})}\\
&\geq \mybox{G_{y_{n+p},\,y_{n+p+1},\,y_{n+p+2}}(\frac{t}{2^{p+1}})\ast G_{y_{n+p+1},\,y_{n+p+2},\,y_{n+p+3}}(\frac{t}{2^{p+2}})\ast\cdots\ast G_{y_{n+q-2},\,y_{n+q-1},\,y_{n+q}}(\frac{t}{2^{q-1}}) \ast G_{y_{n+q-2},\,y_{n+q-1},\,y_{n+q}}(\frac{t}{2^{q-1}})\;\; [\; \text{ as } y_n\neq y_{n+1} \text{ for all } n\in \mathbb{N}\; ] }\\
&\geq \mybox{G_{y_{n+p},\,y_{n+p+1},\,y_{n+p+2}}(\frac{t}{2^{p+1}})\ast G_{y_{n+p+1},\,y_{n+p+2},\,y_{n+p+3}}(\frac{t}{2^{p+2}})\ast\cdots\ast G_{y_{n+q-2},\,y_{n+q-1},\,y_{n+q}}(\frac{t}{2^{q-1}})}\\
&\geq \mybox{G_{y_{0},\,y_{1},\,y_{2}}(k^{-(n+p)}\frac{t}{2^{p+1}})\ast G_{y_{n+p+1},\,y_{n+p+2},\,y_{n+p+3}}(k^{-(n+p+1)}\frac{t}{2^{p+2}})\ast\cdots\ast G_{y_{n+q-2},\,y_{n+q-1},\,y_{n+q}}(k^{-(n+q-2)}\frac{t}{2^{q-1}})}\\
&\geq \mybox{G_{y_{0},\,y_{1},\,y_{2}}(2^{n-1}t)\ast G_{y_{n+p+1},\,y_{n+p+2},\,y_{n+p+3}}(2^{n-1}t)\ast\cdots\ast
 G_{y_{n+q-2},\,y_{n+q-1},\,y_{n+q}}(2^{n-1}t),\;\; \text{ as } k\in(0, \frac{1}{2}]}\\
&\geq G_{y_{0},\,y_{1},\,y_{2}}(2^{n-1}t).
\end{align*}
So, from $(v)$ we get
\begin{align*}
&G_{y_n,\,y_{n+p},\,y_{n+q}}(t)\\
&\geq \mybox{G_{y_n,\,y_{n+1},\,y_{n+1}}(\frac{t}{2})\ast G_{y_{n+1},\,y_{n+2},\,y_{n+2}}(\frac{t}{2^2})\ast \cdots \ast G_{y_{n+p-1},\,y_{n+p},\,y_{n+p}}(\frac{t}{2^p})\ast G_{y_{0},\,y_{1},\,y_{2}}(2^{n-1}t)}\\
&\geq \mybox{G_{y_n,\,y_{n+1},\,y_{n+2}}(\frac{t}{2})\ast G_{y_{n+1},\,y_{n+2},\,y_{n+3}}(\frac{t}{2^2})\ast \cdots \ast G_{y_{n+p-1},\,y_{n+p},\,y_{n+p+1}}(\frac{t}{2^p})\ast G_{y_{0},\,y_{1},\,y_{2}}(2^{n-1}t)}\\
&\geq \mybox{G_{y_0,\,y_{1},\,y_{2}}(k^{-n}\frac{t}{2})\ast G_{y_{0},\,y_{1},\,y_{2}}(k^{-(n+1)}\frac{t}{2^2})\ast \cdots \ast G_{y_{0},\,y_{1},\,y_{2}}(k^{-(n+p-1)}\frac{t}{2^p})\ast G_{y_{0},\,y_{1},\,y_{2}}(2^{n-1}t)}\\
&\geq \mybox{G_{y_0,\,y_{1},\,y_{2}}(2^{n}\frac{t}{2})\ast G_{y_{0},\,y_{1},\,y_{2}}(2^{(n+1)}\frac{t}{2^2})\ast \cdots \ast G_{y_{0},\,y_{1},\,y_{2}}(2^{(n+p-1)}\frac{t}{2^p})\ast G_{y_{0},\,y_{1},\,y_{2}}(2^{n-1}t),\;\; \text{ as } k\in(0, \frac{1}{2}]}\\
&\geq \mybox{G_{y_0,\,y_{1},\,y_{2}}(2^{(n-1)}t)\ast G_{y_{0},\,y_{1},\,y_{2}}(2^{(n-1)}t)\ast \cdots \ast G_{y_{0},\,y_{1},\,y_{2}}(2^{(n-1)}t)\ast G_{y_{0},\,y_{1},\,y_{2}}(2^{n-1}t)}\\
&\geq G_{y_{0},\,y_{1},\,y_{2}}(2^{n-1}t).
\end{align*}
So, taking $n\rightarrow\infty$, we get $G_{y_n,\,y_{n+p},\,y_{n+q}}(t)\rightarrow 1$ for all $t>0$. So $\{y_n\}$ is a Cauchy sequence. Since the space is complete, the sequence $\{y_n\}$ converges to some point $z\in X$. So, the subsequences $Ax_{3n},\, Bx_{3n+1}, Cx_{3n+2},\, Dx_{3n+2},\, Sx_{3n+3},\, Tx_{3n+1}$ of $\{y_n\}$ converge to the point $z$ of $X$.\\
\textbf{Step I: }Since $Ax_{3n}\rightarrow z$,  $Sx_{3n}\rightarrow z$, $[A, S]$ is compatible and $S$ is continuous, by the Theorem \ref{th2}, we have $ASx_{3n}\rightarrow Sz$ and $SSx_{3n}\rightarrow Sz$. So,
\begin{align*}
\begin{aligned}
G_{ASx_{3n},\, Bx_{3n+1},\, Cx_{3n+2}}(kt)
&\geq G_{SSx_{3n},\, Tx_{3n+1},\, Dx_{3n+2}}(t)*G_{SSx_{3n},\, ASx_{3n},\, Dx_{3n+2}}(t)*\\
&\smallskip\hspace{.5 cm}G_{ASx_{3n},\, Bx_{3n+1},\, Cx_{3n+2}}(t)
*G_{Tx_{3n+1},\, Bx_{3n+1},\, Cx_{3n+2}}(t)*\\
&\smallskip\hspace{.5 cm}G_{SSx_{3n},\, Cx_{3n+2},\, Dx_{3n+2}}(2t)\\
&\geq G_{SSx_{3n},\, Tx_{3n+1},\, Dx_{3n+2}}(t)*G_{SSx_{3n},\, ASx_{3n},\, Dx_{3n+2}}(t)*\\
&\smallskip\hspace{.5 cm}G_{ASx_{3n},\, Bx_{3n+1},\, Cx_{3n+2}}(t)*G_{Tx_{3n+1},\, Bx_{3n+1},\, Cx_{3n+2}}(t)*\\
&\smallskip\hspace{.5 cm}G_{SSx_{3n},\, Cx_{3n+2},\, Dx_{3n+2}}(t).
\end{aligned}
\end{align*}
Taking $n\rightarrow\infty$  we get from above
\begin{align*}
\begin{aligned}
G_{Sz,\, z,\, z}(kt)
 &\geq  G_{Sz,\, z,\, z}(t)*G_{Sz,\, Sz,\, z}(t)*G_{Sz,\, z,\, z}(t)*G_{z,\, z,\, z}(t)*G_{Sz,\, z,\, z}(t)\\ &\smallskip\hspace{7cm}\text{ by the Theorem }\ref{th3} ]\\
&\geq \mybox{G_{Sz,\, z,\, z}(t)*G_{Sz,\, Sz,\, z}(t)*G_{Sz,\, z,\, z}(t)*1\;\; \text{ as } G_{z,\, z,\, z}(t)=1}\\
&=\mybox{ G_{Sz,\, z,\, z}(t)*G_{Sz,\, Sz,\, z}(t)*G_{Sz,\, z,\, z}(t)}\\
&\geq\; G_{Sz,\, z,\, z}(t)*G_{Sz,\, Sz,\, z}(t)\\
&\geq\; G_{Sz,\, z,\, z}(t)*G_{Sz,\, z,\, z}(\frac{t}{2})*G_{z,\, Sz,\, z}(\frac{t}{2})\\
&=\; G_{Sz,\, z,\, z}(t)*G_{Sz,\, z,\, z}(\frac{t}{2})*G_{Sz,\, z,\, z}(\frac{t}{2})\\
&\geq G_{Sz,\, z,\, z}(\frac{t}{2}).
\end{aligned}
\end{align*}
Therefore $G_{Sz,\, z,\, z}(kt)\geq\; G_{Sz,\, z,\, z}(\frac{t}{2})$ for all $k\in(0, \frac{1}{2}]$  and for all $t>0.$ \\
 But $G_{Sz,\, z,\, z}(kt)\leq\; G_{Sz,\, z,\, z}(\frac{t}{2})$ for all $k\in(0, \frac{1}{2}]$  and for all $t>0.$\\
So, $G_{Sz,\, z,\, z}(kt)=\; G_{Sz,\, z,\, z}(\frac{t}{2})$ for all $k\in(0, \frac{1}{2}]$ and for all $t>0,$ which implies that\\
 $G_{Sz,\, z,\, z}(t)=1$  for all $t>0.$  Thus  $Sz=z.$
\begin{align*}
\begin{aligned}
&\text{\;\;\textbf{Step II:} }G_{Az,\, Bx_{3n+1},\, Cx_{3n+2}}(kt)&\geq \;& G_{Sz,\, Tx_{3n+1},\, Dx_{3n+2}}(t)*G_{Sz,\, Az,\, Dx_{3n+2}}(t)\\
& &&*G_{Az,\, Bx_{3n+1},\, Cx_{3n+2}}(t)
*G_{Tx_{3n+1},\, Bx_{3n+1},\, Cx_{3n+2}}(t)\\
&&&*G_{Sz,\, Cx_{3n+2},\, Dx_{3n+2}}(t).
\end{aligned}
\end{align*}
Taking $ n\rightarrow\infty $ we get
\begin{align*}
\begin{aligned}
G_{Az,\, z,\, z}(kt) &\geq \; G_{Sz,\, z,\, z}(t)*G_{Sz,\, Az,\, z}(t)*G_{Az,\, z,\, z}(t)*G_{z,\, z,\, z}(t)*G_{Sz,\, z,\, z}(t)\\
&\geq\; G_{z,\, z,\, z}(t)*G_{z,\, Az,\, z}(t)*G_{Az,\, z,\, z}(t)*G_{z,\, z,\, z}(t)*G_{z,\, z,\, z}(t) \;\; [ \text{ as } Sz=z ]\\
&=\; G_{Az,\, z,\, z}(t)*G_{Az,\, z,\, z}(t) \text{ as } G_{z,\, z,\, z}(t)=1 \\
&\geq\; G_{Az,\, z,\, z}(t).
\end{aligned}
\end{align*}
Therefore,\; $G_{Az,\, z,\, z}(kt)\geq\; G_{Az,\, z,\, z}(t)$  for all $k\in(0, \frac{1}{2}]$  and for all $t>0.$ \\
So, $G_{Az,\, z,\, z}(t)=1$  for all $t>0$. Thus  $Az=z.$
 
\textbf{Step III:} Since $Bx_{3n+1}\rightarrow z$,  $Tx_{3n+1}\rightarrow z$,  $[B, T]$ is compatible and $T$ is continuous, by the Theorem \ref{th2} we get $BTx_{3n+1}\rightarrow Tz$ and $TTx_{3n+1}\rightarrow Tz$. So,\\
\begin{align*}
\begin{aligned}
G_{Ax_{3n},\, BTx_{3n+1},\, Cx_{3n+2}}(kt)&\geq \; G_{Sx_{3n},\, TTx_{3n+1},\, Dx_{3n+2}}(t)*G_{Sx_{3n},\, Ax_{3n},\, Dx_{3n+2}}(t)*\\
&\smallskip\hspace{.5 cm}G_{Ax_{3n},\, BTx_{3n+1},\, Cx_{3n+2}}(t)*G_{TTx_{3n+1},\, BTx_{3n+1},\, Cx_{3n+2}}(t)*\\
&\smallskip\hspace{.5 cm}G_{Sx_{3n},\, Cx_{3n+2},\, Dx_{3n+2}}(2t)\\
&\geq G_{Sx_{3n},\, TTx_{3n+1},\, Dx_{3n+2}}(t)*G_{Sx_{3n},\, Ax_{3n},\, Dx_{3n+2}}(t)*\\
&\smallskip\hspace{.5 cm}G_{Ax_{3n},\, BTx_{3n+1},\, Cx_{3n+2}}(t)*G_{TTx_{3n+1},\, BTx_{3n+1},\, Cx_{3n+2}}(t)*\\
&\smallskip\hspace{.5 cm}G_{Sx_{3n},\, Cx_{3n+2},\, Dx_{3n+2}}(t).
\end{aligned}
\end{align*}
Taking $n\rightarrow\infty$ we get
\begin{align*}
\begin{aligned}
 G_{z,\, Tz,\, z}(kt) &\geq G_{z,\, Tz,\, z}(t)*G_{z,\, z,\, z}(t)*G_{z,\, Tz,\, z}(t)*G_{Tz,\, Tz,\, z}(t)*G_{z,\, z,\, z}(t)\\
&\geq G_{z,\, z,\, z}(t)*G_{Tz,\, z,\, z}(t)*G_{Tz,\, Tz,\, z}(t)\\
&=\; 1*G_{Tz,\, z,\, z}(t)*G_{Tz,\, Tz,\, z}(t)\\
&=\; G_{Tz,\, z,\, z}(t)*G_{Tz,\, Tz,\, z}(t)\\
&\geq\; G_{Tz,\, z,\, z}(t)*G_{Tz,\, z,\, z}(\frac{t}{2})*G_{z,\, Tz,\, z}(\frac{t}{2}).
\end{aligned}
\end{align*}
Therefore, $G_{Tz,\, z,\, z}(kt)\geq G_{Tz,\, z,\, z}(\frac{t}{2})$  for all $k\in(0, \frac{1}{2}]$  and for all  $t>0$. \\
So, $G_{Tz,\, z,\, z}(t)=1 $ for all $t>0$.  Thus $ Tz=z.$
\begin{align*}
\text{\;\;\textbf{Step IV:} }G_{Ax_{3n},\, Bz,\, Cx_{3n+2}}(kt)&\geq G_{Sx_{3n},\, Tz,\, Dx_{3n+2}}(t)*G_{Sx_{3n},\, Ax_{3n},\, Dx_{3n+2}}(t)*\\
&\smallskip\hspace{.5 cm}G_{Ax_{3n},\, Bz,\, Cx_{3n+2}}(t)*G_{Tz,\, Bz,\, Cx_{3n+2}}(t)*\\
&\smallskip\hspace{.5 cm}G_{Sx_{3n},\, Cx_{3n+2},\, Dx_{3n+2}}(t).
\end{align*}
Taking $n\rightarrow\infty$ we get
\begin{align*}
 G_{z,\, Bz,\, z}(kt) &\geq  G_{z,\, Tz,\, z}(t)*G_{z,\, z,\, z}(t)*G_{z,\, Bz,\, z}(t)*G_{Tz,\, Bz,\, z}(t)*G_{z,\, z,\, z}(t)\\
& \geq G_{z,\, z,\, z}(t)*G_{z,\, z,\, z}(t)*G_{z,\, Bz,\, z}(t)*G_{z,\, Bz,\, z}(t)*G_{z,\, z,\, z}(t) \;\; [ \text{ as } Tz=z ]\\
&\geq\; G_{z,\, z,\, z}(t)*G_{Bz,\, z,\, z}(t)\\
&=\; 1*G_{Bz,\, z,\, z}(t)\\
&=\; G_{Bz,\, z,\, z}(t).
\end{align*}
Therefore $G_{Bz,\, z,\, z}(kt)\geq\; G_{Bz,\, z,\, z}(t)$  for all $ k\in(0, \frac{1}{2}]$  and for all $ t>0.$ \\
So, $G_{Bz,\, z,\, z}(t)=1$ for all  $t>0$. Thus $ Bz=z.$

\textbf{Step V:} Since $Cx_{3n+2}\rightarrow z$,  $Dx_{3n+2}\rightarrow z$, $[C, D]$ is compatible and $D$ is continuous, $CDx_{3n+2}\rightarrow Dz$ and $DDx_{3n+2}\rightarrow Dz$. So,
\begin{align*}
G_{Ax_{3n},\, Bx_{3n+1},\, CDx_{3n+2}}(kt)&\geq G_{Sx_{3n},\, Tx_{3n+1},\, DDx_{3n+2}}(t)*G_{Sx_{3n},\, Ax_{3n},\, DDx_{3n+2}}(t)*\\
&\smallskip\hspace{.5 cm}G_{Ax_{3n},\, Bx_{3n+1},\, CDx_{3n+2}}(t)*G_{Tx_{3n+1},\, Bx_{3n+1},\, CDx_{3n+2}}(t)*\\
&\smallskip\hspace{.5 cm}G_{Sx_{3n},\, CDx_{3n+2},\, DDx_{3n+2}}(t).
\end{align*}
\begin{align*}
G_{z,\, z,\, Dz}(kt) &\geq  G_{z,\, z,\, Dz}(t)*G_{z,\, z,\, Dz}(t)*G_{z,\, z,\, Dz}(t)*G_{z,\, z,\, Dz}(t)*G_{z,\, Dz,\, Dz}(t)\\
&\geq G_{Dz,\, z,\, z}(t)*G_{Dz,\, Dz,\, z}(t)\\
&\geq G_{Dz,\, z,\, z}(t)*G_{Dz,\, z,\, z}(\frac{t}{2})*G_{z,\, Dz,\, z}(\frac{t}{2}).
\end{align*}
Therefore $G_{Dz,\, z,\, z}(kt)\geq G_{Dz,\, z,\, z}(\frac{t}{2})$ for all  $k\in(0, \frac{1}{2}]$  and for all  $t>0.$ \\
So, $G_{Dz,\, z,\, z}(t)=1$ for all $t>0.$ Thus $Dz=z.$
\begin{align*}
\text{ \textbf{Step VI:} } G_{Ax_{3n},\, Bx_{3n+1},\, Cz}(kt)&\geq  G_{Sx_{3n},\, Tx_{3n+1},\, Dz}(t)*G_{Sx_{3n},\, Ax_{3n},\, Dz}(t)*\\
&\smallskip\hspace{.5 cm} G_{Ax_{3n},\, Bx_{3n+1},\, Cz}(t)*G_{Tx_{3n+1},\, Bx_{3n+1},\, Cz}(t)*\\
&\smallskip\hspace{.5 cm}G_{Sx_{3n},\, Cz,\, Dz}(t).
\end{align*}
Taking $n\rightarrow\infty$ we get
\begin{align*}
 G_{z,\, z,\, Cz}(kt) &\geq  G_{z,\, z,\, Dz}(t)*G_{z,\, z,\, Dz}(t)*G_{z,\, z,\, Cz}(t)*G_{z,\, z,\, Cz}(t)*G_{z,\, Cz,\, Dz}(t)\\
&\geq G_{z,\, z,\, z}(t)*G_{z,\, z,\, z}(t)*G_{z,\, z,\, Cz}(t)*G_{z,\, z,\, Cz}(t)*G_{z,\, Cz,\, z}(t) \;\; [ \text{ as } Dz=z ]\\
&\geq G_{z,\, z,\, z}(t)*G_{Cz,\, z,\, z}(t)\\
&= 1*G_{Cz,\, z,\, z}(t)\\
&= G_{Cz,\, z,\, z}(t).
\end{align*}
Therefore $G_{Cz,\, z,\, z}(kt)\geq G_{Cz,\, z,\, z}(t)$  for all  $k\in(0, \frac{1}{2}]$ and for all  $t>0.$ \\
So,  $G_{Cz,\, z,\, z}(t)=1$  for all $t>0.$  Thus  $Cz=z.$\\ 
So, we get $Az=Bz=Cz=Dz=Sz=Tz=z$. Therefore $z$ is a common fixed point of $A,\, B,\, C,\, D,\, S$ and $T$. For Uniqueness, let $w$ be another fixed point of these functions. Then we have
\begin{align*}
&G_{Az,\, Bz,\, Cw}(kt)&\geq \;& G_{Sz,\, Tz,\, Dw}(t)*G_{Sz,\, Az,\, Dw}(t)*G_{Az,\, Bz,\, Cw}(t)\\
& &&*G_{Tz,\, Bz,\, Cw}(t)*G_{Sz,\, Cw,\, Dw}(t)\\
&\text{or},\; G_{z,\, z,\, w}(kt) &\geq \;& G_{z,\, z,\, w}(t)*G_{z,\, z,\, w}(t)*G_{z,\, z,\, w}(t)*G_{z,\, z,\, w}(t)*G_{z,\, w,\, w}(t)\\
&&\geq\;& G_{w,\, z,\, z}(t)*G_{w,\, w,\, z}(t)\\
&&\geq\;& G_{w,\, z,\, z}(t)*G_{w,\, z,\, z}(\frac{t}{2})*G_{z,\, w,\, z}(\frac{t}{2})\\
&&\geq\;& G_{w,\, z,\, z}(\frac{t}{2}).
\end{align*}
Therefore $G_{w,\, z,\, z}(kt)\geq\; G_{w,\, z,\, z}(\frac{t}{2})$ for all $k\in(0, \frac{1}{2}]$  and for all $t>0.$ \\
So, $G_{w,\, z,\, z}(t)=1$ for all $t>0.$ Thus $w=z.$
Therefore $z$ is a unique fixed point of $A,\, B,\, C,\, D,\, S$ and $T$.

\end{document}